\documentclass[11pt,a4paper]{article}
\usepackage[latin1]{inputenc}
\usepackage{graphicx}
\usepackage[mathscr]{euscript}
\usepackage{a4}
\usepackage{epsfig}
\usepackage{stmaryrd,amsfonts,amsmath,amssymb}
\usepackage{ifthen}
\usepackage{hyperref}

\newcommand{\origforall}{}    \let\origforall=\forall
\renewcommand{\forall}{\,\origforall\,}
\newcommand{\origexists}{}    \let\origexists=\exists
\renewcommand{\exists}{\,\origexists\,}
\newcommand{\origin}{}        \let\origin=\in
\renewcommand{\in}{{\,\origin\,}}
\newcommand{\origwedge}{}     \let\origwedge=\wedge
\renewcommand{\wedge}{{\quad\origwedge\quad}}
\newcommand{\origvee}{}       \let\origvee=\vee
\renewcommand{\vee}{{\quad\origvee\quad}}
\newcommand{\origfrac}{}      \let\origfrac=\frac
\renewcommand{\frac}[2]{{\;\origfrac{#1}{#2}\;}}
    
\renewcommand{\subset}{\subseteq}

\newcommand{\beq}{\begin{eqnarray*}}
\newcommand{\eeq}{\end{eqnarray*}}

\newcommand{\hence}{\curvearrowright}
\newcommand{\R}{\mathbb{R}}

\newcommand{\N}{\mathbb{N}}
\newcommand{\Z}{\mathbb{Z}}

\newcommand{\n}{\textnormal}

\renewcommand{\S}{\mathfrak{S}}
\DeclareMathOperator{\Lip}{{\n{Lip}}}
\newcommand{\smax}{\origvee}

\DeclareMathOperator{\lcm}{\n{lcm}}

\DeclareMathOperator{\im}{\n{im}}
\DeclareMathOperator{\Cay}{\n{Cay}}
\DeclareMathOperator{\diam}{\n{diam}}
\DeclareMathOperator{\Isom}{\n{Isom}}
\DeclareMathOperator{\eIsom}{\n{-Isom}}
\DeclareMathOperator{\URIsom}{\n{URIsom}}
\DeclareMathOperator{\UQIsom}{\n{UQIsom}}
\DeclareMathOperator{\eIden}{\n{-Iden}}
\DeclareMathOperator{\Iden}{\n{Iden}}
\DeclareMathOperator{\id}{\n{id}}
\DeclareMathOperator{\GL}{\n{GL}}
\newcommand{\fg}[1]{C_{#1}}

\newcommand{\imetric}{me\-tric~}  

\newcommand{\ig}[1]{}

\newcommand{\esquare}{\hfill\ensuremath{\square}}

\newcounter{claim}
\newcounter{titleclaim}
\newenvironment{Proof}{{\textbf{Proof}$\;$}}{\esquare\newline}
\newenvironment{Non-Proof}{{\textbf{Non-Proof}$\;$}}{\hfill\ensuremath{\times}\newline}
\newenvironment{History}{{\textbf{Historical Remark}$\;$}}{}

\newenvironment{Question}{{\textbf{Question}$\;$}\begin{itshape}}{\end{itshape}}

\newtheorem{LemmaBase}[claim]{Lemma}
\newenvironment{Lemma}[1][-1]{\ifthenelse{\equal{#1}{-1}}
    {\begin{LemmaBase}}{\begin{LemmaBase}[#1]}
    \dotfill \, {\em \bf \arabic{claim}} \quad\\{} }
    {\addtocounter{titleclaim}{1}\end{LemmaBase}}

\newtheorem{CorollaryBase}[claim]{Corollary}
\newenvironment{Corollary}[1][-1]{\ifthenelse{\equal{#1}{-1}}
    {\begin{CorollaryBase}}{\begin{CorollaryBase}[#1]}
    \dotfill \, {\em \bf \arabic{claim}} \quad\\{} }
    {\addtocounter{titleclaim}{1}\end{CorollaryBase}}

\newtheorem{PropositionBase}[claim]{Proposition}
\newenvironment{Proposition}[1][-1]{\ifthenelse{\equal{#1}{-1}}
    {\begin{PropositionBase}}{\begin{PropositionBase}[#1]}
    \dotfill \, {\em \bf \arabic{claim}} \quad\\{} }
    {\addtocounter{titleclaim}{1}\end{PropositionBase}}

\newtheorem{ConjectureBase}[claim]{Conjecture}

\newtheorem{ExampleBase}[claim]{Example}
\newenvironment{Example}[1][-1]{\ifthenelse{\equal{#1}{-1}}
    {\begin{ExampleBase}}{\begin{ExampleBase}[#1]}
    \dotfill \, {\em \bf \arabic{claim}} \quad\\{} }
    {\addtocounter{titleclaim}{1}\end{ExampleBase}}

\newtheorem{DefinitionBase}[claim]{Definition}
\newenvironment{Definition}[1][-1]{\ifthenelse{\equal{#1}{-1}}
    {\begin{DefinitionBase}}{\begin{DefinitionBase}[#1]}
    \dotfill \, {\em \bf \arabic{claim}} \quad\\{} }
    {\addtocounter{titleclaim}{1}\end{DefinitionBase}}

\newtheorem{TheoremBase}[claim]{Theorem}
\newenvironment{Theorem}[1][-1]{\ifthenelse{\equal{#1}{-1}}
    {\begin{TheoremBase}}{\begin{TheoremBase}[#1]}
    \hrulefill \, {\em \bf \arabic{claim}} \quad\\{} }
    {\addtocounter{titleclaim}{1}\end{TheoremBase}}

\begin{document}

\title{Shared Rough and Quasi-Isometries of Groups}
\author{Andreas Lochmann}
\maketitle

\abstract

We present a variation of quasi-isometry to approach the problem of defining a
geometric notion equivalent to commensurability. In short, this variation can be
summarized as ``quasi-isometry with uniform parameters for a large enough
family of generating systems''. Two similar notions (using isometries and rough
isometries instead, respectively) are presented alongside.

This article is based mainly on a chapter of the author's doctoral thesis
(\cite{Lochmann_dissertation}).

\section{Basic Notions}

\subsection{Rough and Quasi-Isometries}

A well written introduction to coarse geometry is the book by Burago, Burago
and Ivanov (\cite{Burago}). In the following, let $X$ and $Y$ be \imetric spaces.

\begin{Definition} \label{DEF_near}
Two (set theoretic) mappings $\alpha,\beta: X\rightarrow Y$ are {\em
$\epsilon$-near} to each other, $\epsilon \geq 0$, if $d_Y(\alpha x, \beta x)
\leq \epsilon \forall x\in X$. (We drop brackets where feasible.)

A (set theoretic) mapping $\alpha: X\rightarrow Y$ is {\em
$\epsilon$-surjective}, $\epsilon \geq 0$, if for each $y\in Y$ there is
$x\in X$ such that $d_Y(\alpha x, y) \leq \epsilon$.
\end{Definition}

\begin{Definition}\label{DEF_quasi_isometry}
A (not neccessarily continuous) map $\eta: X\rightarrow Y$ is called a {\em
$(\lambda,\,\epsilon)$-quasi-isometric embedding}, $\epsilon,\,\lambda\, \geq\,
0$ (which shall always imply $\lambda,\,\epsilon\,\in\,\R$), if
\beq
\lambda^{-1}\,d_X(x,\,x') \,-\,\epsilon \;\;\leq\;\; d_Y(\eta x,\, \eta x')
\;\;\leq\;\; \lambda\,d_X(x,\,x') \,+\,\epsilon
\eeq
for all $x, x'\in X$. 

A pair $\eta: X\rightarrow Y$, $\eta': Y\rightarrow X$ of
$(\lambda,\,\epsilon)$-quasi-isometric embeddings is called a {\em
$(\lambda,\,\epsilon)$-quasi-isometry} if $\eta\circ\eta'$ and $\eta'\circ
\eta$ are $\epsilon$-near the identities on $Y$ and $X$, respectively. When we
speak of a ``quasi-isometry $\eta:X\rightarrow Y$'' a corresponding map $\eta'$
shall always be implied.

$X$ and $Y$ are called quasi-isometric, if there is a
quasi-isometry between them. 
\end{Definition}

\begin{Definition}\label{DEF_rough_isometry}
A $(1,\,\epsilon)$-quasi-isometric embedding, $\epsilon \geq 0$, is called an
{\em $\epsilon$-isometric embedding}. It fulfills
\beq
|d_X(x,x') - d_Y(\eta x, \eta x')| &\leq& \epsilon
\eeq
for all $x, x'\in X$. 

An {\em $\epsilon$-isometry} is a $(1,\,\epsilon)$-quasi-isometry.
The map $\eta:X\rightarrow Y$ is called a {\em rough isometry} if there is some
$\epsilon \geq 0$ such that $\eta$ is an $\epsilon$-isometry.

$X$ and $Y$ are called $\epsilon$-isometric [roughly isometric], if there is
an $\epsilon$-isometry [any $\epsilon$-isometry] between them. 
\end{Definition}

\begin{History}
It is difficult to attribute the concept of rough isometry to a single person,
as it was always present in the notion of quasi-isometry, which itself was an
obvious generalization of what was then called pseudo-isometry by Mostow in his
1973-paper about rigidity (see \cite{Mostow}, \cite{Gromov_Hyperbolic},
\cite{Kanai}). Recent developments about the stability of rough isometries can
be found in \cite{Rassias_Isometries}.
\end{History}

\ig{
\begin{Definition}
A mapping $f: X\rightarrow Y$ is called {\em $(K,\, \epsilon)$-Lipschitz}
(i.e.\ ``$K$-Lipschitz map on $\epsilon$-scale'' in \cite{Gromov_Asymptotic}),
$\epsilon,\, K\,\geq\, 0$, iff 
\beq
d_Y(f(x),\, f(y)) &\leq& K\cdot d_X(x,\,y) \;+\; \epsilon \qquad \forall
x,\,y\,\in\, X.
\eeq
If $\epsilon\,=\,0$, $f$ is {\em $K$-Lipschitz (continuous)}. Define
$\Lip_{K,\,\epsilon}(X,\,Y)$ to be the set of all $(K,\,\epsilon)$-Lipschitz
functions $X\rightarrow Y$, and $\Lip_{K,\,\epsilon} X :=
\Lip_{K,\,\epsilon}(X,\,[0,\,\infty])$, $\Lip X := \Lip_{1,0}(X)$.
\end{Definition}

If nothing else is said, $[0,\,\infty] = \R_{\geq 0}\cup\{\infty\}$ is the default
codomain for a Lipschitz function.

Assume $f$ to be a $(K,\,\epsilon)$-Lipschitz function on $X$ and $f(x)\, =\,
\infty$ for some $x\in X$. Then clearly $f(y)\, =\, \infty$ for all $y$ in
finite distance to $x$. Thus, if $X$ is a true metric space, we have $\Lip X 
= \Lip(X,\,\R_{\geq 0}) \cup \{\infty\}$.
}

\subsection{Commensurability}

When one speaks about the coarse geometry of finitely generated groups, one
generally means quasi-isometries of Cayley graphs. While a single infinite group
gives rise to an infinite number of non-isomorphic Cayley graphs,
quasi-isometries do not depend on the generating system of the group, and hence
the quasi-isometry class of a group is well-defined, and an important invariant.
It encompasses the idea of two groups being {\em approximately isomorphic}, but
quasi-isometries are not the only way to do this. Particularly the pure
group-theoretic notion of commensurability rivals the quasi-isometry, and their
interplay is still an interesting research problem. In the following, we use the
definitions given in \cite{dlHarpe}.

\begin{Definition} \label{DEF_commensurable}
Let $G$ and $H$ be groups. $G$ and $H$ are {\em commensurable} when there exist
subgroups $G' \leq G$ and $H' \leq H$ of finite index, such that $G'$ and $H'$ are
isomorphic as group.

$G$ and $H$ are {\em commensurable up to finite kernels} if there exists
a finite sequence of groups $\Gamma_1,\, \ldots,\, \Gamma_N$ and homomorphisms
$h_0,\,\ldots,\, h_N$ 
\beq
G \;\stackrel{h_0}{\longrightarrow}\;
\Gamma_1 \;\stackrel{h_1}{\longleftarrow}\;
\Gamma_2 \;\stackrel{h_2}{\longrightarrow}\;
\Gamma_3 \;\stackrel{h_3}{\longleftarrow}\;
\;\ldots\; \;\stackrel{h_{N-1}}{\longrightarrow}\;
\Gamma_N \;\stackrel{h_N}{\longleftarrow}\; H
\eeq
with finite kernels and images of finite index.
\end{Definition}

One easily sees that commensurability always implies commensurability up to
finite kernels, which in turn always implies quasi-isometry, given that both
groups are finitely generated. We quote without proof the following Proposition
from \cite{dlHarpe}, IV.28.

\begin{Proposition}
Two residually finite groups are commensurable if and only if they are
commensurable up to finite kernels.
\end{Proposition}

\begin{Proposition} \label{PRO_homomorphism_and_qi}
Let $G$ and $H$ be f.g.\ groups, and $\eta:\,G\rightarrow H$ a homomorphism and
quasi-isometry. Then $G$ and $H$ are commensurable up to finite kernels.
\end{Proposition}
\begin{Proof}
The kernel of $\phi$ is finite, because it is the preimage of a
finite subset of $H$. And the image $\phi(G)$ is a subgroup of
$H$ of finite index: $\phi(G)$ is $\epsilon$-dense in $H$. Let $B$ be the
$\epsilon$-ball around the identity in $H$, then each element $h\,\in\,H$ can be
written as $b\cdot \phi(g)$ for some $b\,\in\,B$ and $g\,\in\,G$. With this, the
number of cosets of $G/\phi(g)$ can be at most as large as $\# B$, and in
particular, it is finite.
\end{Proof}

There is a multitude of cases in which quasi-isometry implies commensurability
(for example f.g.\ abelian groups, certain types of Baumslag-Solitar groups,
abelian-by-cyclic groups in \cite{Farb_Mosher}) but also a plenty supply of
counter-examples (e.g.\ Lamplighter groups, or $\Z^2\,\rtimes_A\, \Z$ with certain
choices for $A\,\in\, \GL(2,\,\Z)$).

One would think that strengthening the geometric equivalence to rough isometry
could be enough to imply commensurability---this, however, is wrong: In
contrast to quasi-isometry, rough isometry is not a canonical notion for
groups, as it depends on the chosen generating set (see Section
\ref{SEC_exponential_growth_rate}). The very next idea would be to ask for the
{\it existence} of generating systems of the two groups, such that they are
roughly isometric. This doesn't lead to commensurability as well, the
lamplighter groups (see \cite{dlHarpe} IV.44) present counter-examples.

\ig{
A related question is whether two isomorphic Cayleygraphs imply that their groups
are isomorphic. This is in general not the case (as long as only one Cayleygraph
per group is considered), and is already in the finite case a rich source for
problems, see \cite{Li}.
}

\ig{
We are not yet able to categorize all of these notions and give proofs or
counter-examples to their mutual equivalences. However, we will take a closer
look at two of these definitions. Our methods involve the analysis of generating
systems and groups of what we call ``shared isometries'' and ``shared rough
isometries'' -- maps which are isometries (respectively rough isometries)
relative to lots of generating systems at once. But before we get there, we
insert a section about a rough isometry invariant, and another section about
the case in abelian groups. These two sections together provide first insights.
}

Our approach will be to require the existence of a sufficiently large class of
quasi-isometries with bounded parameters $(\lambda,\,\epsilon)$ for a whole
family of generating systems.

\subsection{Notation}

Denote with:
\begin{itemize}
\item $\fg{n}$ the cyclic group of order $n$,
\item $\Cay(G,\,S)$ the Cayley graph of the group $G$ relative to the generating
system $S\,\subset\,G$, and
\item $\N^*$ the natural numbers without zero.
\end{itemize}


\section{Exponential Growth Rate} \label{SEC_exponential_growth_rate}

Each quasi-isometry invariant is also a rough isometry invariant. But there
also is a rough isometry invariant, which is not a quasi-isometry invariant;
this can be used to easily deny the existence of certain rough isometries.

\begin{Definition}
Let $G$ be a f.g.\ group, and let $S$ be a finite generating system of $G$. The
{\em exponential growth rate} is
\beq
\omega(G,\, S) &:=& \limsup_{k\rightarrow \infty} \;\sqrt[k]{\#B_{G,S}(k)} \;\;=\;\;
\exp \; \limsup_{k\rightarrow\infty}\, \frac{\,\ln \#B_{G,S}(k)\,}{k}.
\eeq
The {\em minimal growth rate} $\omega(G)$ is the infimum of $\omega(G,S)$ over
all finite generating systems $S$. The group $G$ is of {\em uniformly exponential growth}
if $\omega(G) > 1$.
\end{Definition}

\begin{Proposition}
Let $F_n$ be the free group on $n$ generators. Then $\omega(F_n) = 2\,n-1$.
\end{Proposition}
\begin{Proof}
This is Proposition VII.12 in \cite{dlHarpe}, we summarize the proof here: The
minimal growth rate is attained by any free generating system for $F_n$.
Now let $S$ be any generating system of $F_n$. Let $S'$ be the image of $S$
under abelianization, choose a minimal subset of $S'$ generating a finite index
subset of $\Z^n$. Any preimage of $S'$ is a set of free generators of a
subgroup $H$ of $F_n$, which in turn is isomorphic to $F_n$ and of growth $2\,n-1$.
Hence, $\omega(F_n,\, S) \,\geq\, \omega(H,\, S') \,=\, 2\,n-1$.
\end{Proof}

\begin{Lemma}
Let $G, H$ be f.g.\ groups of uniformly exponential growth, let $S_G$ and $S_H$ be
finite generating systems of $G$ and $H$, and let
\beq
\eta: \quad \Cay(G,\, S_G) &\rightarrow& \Cay(H,\, S_H)
\eeq
be an $\epsilon$-isometry, $\epsilon \geq 0$. Then
$\omega(G,\, S_G) \,=\, \omega(H,\, S_H)$.
\end{Lemma}
\begin{Proof}
By estimating the number of elements in each ball:
\beq
\frac{\#B_G(r)}{\#B_G(\epsilon)}&\leq& \#\eta(B_G(r)) \\
&\leq& \# B_H(r +
\epsilon)\\
&\leq& \# B_H(r) \cdot \#B_H(\epsilon)\\
\Rightarrow \quad\omega(G,\, S_G) &\leq& \omega(H,\, S_H),
\eeq
and vice versa.
\end{Proof}

\begin{Example} \label{EXA_free_groups_roughly_isometric}
Although the free groups $F_2$ and $F_4$ are commensurable, there is a
generating system of $F_2$, such that its Cayley graph is not roughly isometric
to any Cayley graph of $F_4$, as the minimal growth rates differ.

However, there still exist generating systems of $F_2$ and $F_4$ with their
Cayley graphs being roughly isometric: Choose any embedding $\pi$ of $F_4$ into
$F_2$ as subgroup of finite index, let $S_j$ be free generating systems of
$F_j$, $j=2,4$. Choose $S \,:=\, S_2\,\cup\,\pi(S_4)$ as generating system
for $F_2$, then due to the uniqueness of each word in $F_2$ and due to the
corresponding unique length function, $\pi$ is a rough isometry between
$\Cay(F_2,\,S)$ and $\Cay(F_4,\,S_4)$.
\end{Example}

\section{Rough Isometries of Quotients with\\ Finite Kernel}

In this section we deduce the existence of rough isometries between groups and
their quotients of finite index.

\begin{Proposition}\label{PRO_normal_subgroup}
Let $G$ be a f.g.\ group, let $H\,\trianglelefteq\,G$ be a finite normal
subgroup, and set $G'\,:=\,G/H$. Then each finite generating system $S_0$ of $G$ can be enlarged
to a finite generating system $S$ of $G$, such that $\eta: G\rightarrow
G'$, $g\mapsto gH$ is a 1-isometry, where $G'$ is endowed with the word
metric of the projection $S'$ of $S$.
\end{Proposition}
\begin{Proof}
Let $S_0$ be some generating set of $G$, and put $S \,:=\, S_0 \,\cup\, H \setminus
\{e\}$. Define $S' \,:=\, \{sH\,:\, s\,\in\, S\}\setminus \{e\}$ as the non-trivial
cosets of $S$. $S'$ generates $G'$: For each $x\,\in\, G'$ is $x \,=\, gH$ for some
$g\,\in\, G$, present $g$ as $s_1\ldots s_n$ with $s_j\,\in\, S_0$. Then $x =
(s_1H)\cdot\ldots\cdot (s_nH)$. From this we see $d'(eH,\, gH)\,\leq\, d(e,\,g)$
with $d$ the word metric resulting from $S\subset G$ and $d'$ the word metric
for $S'\,\subset\, G'$.

On the other hand, let $g\,\in\,G$ be arbitrary, and let
$gH\,=\,(s_1H)\cdot\ldots\cdot(s_nH)\,\in\, G'$ be a shortest word in $G'$. Then
there is $h\in H$ with $g = s_1\ldots s_n\cdot h$, hence $d(e,\,g)\,\leq\,
d'(eH,\, gH) \,+\, 1$. Finally, let $x\,\in\, G'$ be any coset. Choose any
representative $g$ of this coset, thus $gH\,=\, x$. Then $d(x,\, \eta(g)) \,=\,
d(x,\, gH) \,=\, 0$.
\end{Proof}

Note that in the preceding proof we might have chosen some generating set $S_H$
of $H$ and set $S := S_0\cup S_H$. In this case, the proof would yield an
$\epsilon$-isometry with $\epsilon = \diam \Cay(H, S_H)$ instead.

\begin{Proposition}\label{PRO_larger_generating_set}
Let $S$ be some generating set of the finitely generated group $G$,
$H\,\trianglelefteq\, G$ finite, and $S_H$ a generating set of
$H$. Then the identity $(G,\, d_S) \rightarrow (G,\, d_{S\,\cup\, S_H})$ is an
$\epsilon$-isometry with $\epsilon \,\leq\, \diam_S(H)$.  
\end{Proposition}
\begin{Proof}
Let $d\,:=\, d_S$, $d' \,:=\, d_{S\,\cup\, H}$. We obviously have $d'(e,\,g)
\,\leq\, d_{S\,\cup\, S_H}(e,\,g) \,\leq\, d(e,\,g)$ for all $g\,\in\, G$. Now let $g\,
=\, s_1\,t_1\ldots s_n\, t_n$ be some presentation of $g\,\in\, G$ in generators
$s_j\,\in\, S\,\cup\,\{e\}$ and $t_j\,\in\, H$. As $H$ is normal, we can find
$t'_1$ to $t'_n\,\in\, H$ with $g\,=\,s_1\ldots s_n\,\cdot\, t'_1\ldots t'_n$.
Hence $d(e,\,g) \,\leq\, d'(e,\,g) \,+\, \epsilon$ where $\epsilon$ is the
diameter of $H\,\subset\, G$ in $d_S$. 
\end{Proof}

\section{Shared Isometries}

\begin{Definition}
Consider $\epsilon \geq 0$, and let $G$ be a finitely generated group. Let $\S$
be a family of generating systems. 
Define the {\em $\S$-shared} or simply {\em shared isometry groups and sets}
\beq
(\lambda,\,\epsilon)\eIsom_\S(G) &:=& \{\eta: G \rightarrow G \;|\; \forall S\in \S:\;\eta
\n{ is a } (\lambda,\,\epsilon)\n{-qi.} \n{ rel. to }S\}\\
\epsilon\eIsom_\S(G) &:=& (1,\,\epsilon)\eIsom_\S(G)\\
\Isom_\S(G) &:=& 0\eIsom_\S(G)\\
\UQIsom_\S(G) &:=& \bigcup_{\lambda,\epsilon \,\geq\, 0}
(\lambda,\,\epsilon)\eIsom_\S(G)\\
\URIsom_\S(G) &:=& \bigcup_{\epsilon \,\geq\, 0} \epsilon\eIsom_\S(G)
\eeq
The last ones we call {\em $\S$-uniform} quasi-isometries resp. rough
isometries. We further define
\beq
\epsilon\eIden_\S(G) &:=& \{\eta: G\rightarrow G \;|\; \forall
S\in\S:\; \eta \n{ is $\epsilon$-near the identity}\}\\
\Iden_\S(G) &:=& \bigcup_{\epsilon\geq 0} \epsilon\eIden_\S(G).
\eeq
\end{Definition}

These definitions are similar to the definition of the quasi-isometry group
$\n{QI}$ of a metric space or group (the calculation of $\n{QI}$ is very
difficult in general, see for example \cite{Farb_Mosher}), and we find
composition to be a group structure on $\UQIsom_\S(G)$ and on $\URIsom_\S(G)$
after quotiening out $\Iden_\S(G)$. The difference between the quasi-isometry
group $\n{QI}(G)$ and $\UQIsom_\S(G)/\Iden_\S(G)$ seems to be subtle,
as we just demand $\lambda$ and $\epsilon$ to be uniformly bounded for all word
metrics in $\S$, but this difference can be enormous, if $\S$ is chosen large
enough. On the other hand, if $\S$ comprises only a finite number of generating
systems, $\UQIsom_\S(G)/\Iden_\S(G)$ equals $\n{QI}(G)$, independently of the
exact choice of $\S$. We will begin with the examination of $\UQIsom_\S(G)$ and
$\URIsom_\S(G)$ in Section \ref{SEC_shared_rough}, and now concentrate on the
nearly trivial case of $\Isom_\S(G)$. We start with a simple observation, which
resulted from a discussion with Laurent Bartholdi and Martin Bridson during the
2007 winter school ``Geometric Group Theory'' in G\"ottingen:

\begin{Theorem} \label{THE_bartholdi}
\indent
{\bf (A)} Let $\S = \S_\n{asym}$ be the family of all, possibly asymmetric, finite
generating systems of $G$. Then $\Isom_\S(G)$ is isomorphic to $G$ (using
possibly asymmetric distance functions).

{\bf (B)} Let $G$ be a group with a finite, symmetric generating
system $S_0$ such that the following hold:
\begin{enumerate}
\item There are no $s_1,s_2,s_3\in S_0$ with $s_1s_2 = s_3$. (Minimality; easy to achieve.)
\item There are no $s_1,s_2\in S_0$, $s_1\neq s_2^{\pm 1}$, with $s_1^2s_2^2 = e$.
\item There are no $s_1,s_2\in S_0$, $s_1\neq s_2^{\pm 1}$, with $s_1^{s_2} = s_1^{-1}$.
\item There are no $s_1,s_2\in S_0$, $s_1\neq s_2^{\pm 1}$, with $s_1^{s_2} = s_1$ \\
(In particular, $G$ is not an abelian group.)
\item There are at least two distinct elements in $S_0$, which are not
inverses of each other.
\end{enumerate}
Let $\S = \S_\n{sym}$ be the family of all symmetric finite generating systems
of $G$. Then $\Isom_\S(G)$ is isomorphic to $G$.

{\bf (C)} Let $G$ be a f.g.\ abelian group without 2-torsion, and let $\S =
\S_\n{sym}$ be the family of all symmetric finite generating systems of $G$.
Then $\Isom_\S(G)$ is isomorphic to $G\rtimes\fg{2}$, where $\fg{2}$ acts by
inversion $x\mapsto x^{-1}$.

{\bf (D)} Let $G$ be a f.g.\ group, and $S_0 \in \S = \S_\n{sym}(G)$, such that $S_0$
is minimal, and each element $s\in S$ has order 2 (i.e.\ $s^2 = e$). Then
$\Isom_\S(G) \cong G$.
\end{Theorem}
\begin{Proof}
The proof is based on an idea by L.\ Bartholdi.

{\bf (A)} Consider $\phi\in \Isom_\S(G)$, and $x, s\in G$ arbitrary, $s\neq e$. Let
$\S' := \{S\in \S: s\in S\}$. Then $d_S(x,xs) = 1$ and $d_S(\phi (x), \phi (xs))
= 1$ for each $S\in\S'$, i.e.\ $s_x := \phi(x)^{-1}\cdot \phi(xs) \in S$. Assume
$s_x\neq s$. Then define $S' := (S\setminus \{s_x\}) \cup \{s,\, s^{-1}s_x\}$.
$S'$ is again a generating system and $s_x\notin S'$, as $s\neq s_x$ and $s\neq
e$. Yet, we have $s\in S'$, contradiction. So we conclude $s_x = s$ and
$\phi(xs) = \phi(x) \cdot s$. By induction we find $\phi(x) = \phi(e) \cdot x$,
with $\phi(e)$ arbitrary. On the other hand, each such $\phi$ obviously is in
$\Isom_\S(G)$, and
\beq
G \,\ni\, g &\mapsto& (\phi_g: \; x\,\mapsto\, g\cdot x) \,\in\,\Isom_\S(G)
\eeq
are shared isometries, and $\phi_g \circ \phi_h \,=\, \phi_{gh}$.

{\bf (B)} Let $\phi\in \Isom_\S(G)$, and $x\in G$ arbitrary, $s\in S_0$.
Then $d_{S_0}(x,xs) = 1$ and $d_{S_0}(\phi(x), \phi(xs)) = 1$, i.e.\
$s_x := \phi(x)^{-1}\cdot \phi(xs)\in S_0$. Like in the asymmetric case, using
$\S' := \{S\in \S: s\in S\} \ni S_0$ we find $s_x = s$ or $s_x = s^{-1}$, but
the choice might depend on $x$, and this is the main point differing to the
asymmetric case. Now let $r\in S_0$ be arbitrary, $r\neq s^{\pm 1}$ and $S'_0 :=
S_0 \cup \{sr, (sr)^{-1}\}$. Note that $d_{S_0}(x, xsr) = 2$, as there are no
triangles in $S_0$, but $d_{S'_0}(x, xsr) = 1$. Let $r_y = \phi(y)^{-1}
\cdot \phi(yr)\in S_0$, so we find $\phi(xsr) = \phi(x)\cdot s_x \cdot r_{xs}$.
As $d_{S'_0}(\phi(x), \phi(xsr)) = 1$, we have
\begin{enumerate}
\item $s_x = s$ or $s_x = s^{-1}$,
\item $r_{xs} = r$ or $r_{xs} = r^{-1}$,
\item $s_xr_{xs} \in S'_0$, but $s_xr_{xs}\notin S_0$.
\end{enumerate}
Hence, $s_xr_{xs}$ must be one of the added elements $sr$ or $(sr)^{-1} =
r^{-1}s^{-1}$. We find eight cases:
\begin{enumerate}
\item $s_x = s,\; r_{xs} = r,\; s_xr_{xs} = sr$
\item $s_x = s^{-1},\; r_{xs} = r,\; s_xr_{xs} = sr
\quad\Rightarrow\quad s^2 = e \quad\hence\quad \n{case (1)}$
\item $s_x = s,\; r_{xs} = r^{-1},\; s_xr_{xs} = sr
\quad\Rightarrow\quad r^2 = e \quad\hence\quad \n{case (1)}$
\item $s_x = s^{-1},\; r_{xs} = r^{-1},\; s_xr_{xs} = sr
\quad\Rightarrow\quad s^2r^2 = e$
\item $s_x = s,\; r_{xs} = r,\; s_xr_{xs} = r^{-1}s^{-1}
\quad\Rightarrow\quad (sr)^2 = e \quad\hence\quad \n{case (1)}$
\item $s_x = s^{-1},\; r_{xs} = r,\; s_xr_{xs} = r^{-1}s^{-1}
\quad\Rightarrow\quad r^s = r^{-1}$
\item $s_x = s,\; r_{xs} = r^{-1},\; s_xr_{xs} = r^{-1}s^{-1}
\quad\Rightarrow\quad s^r = s^{-1}$
\item $s_x = s^{-1},\; r_{xs} = r^{-1},\; s_xr_{xs} = r^{-1}s^{-1}
\quad\Rightarrow\quad r^s = r$
\end{enumerate}
Cases (2), (3) and (5) directly lead to case (1) after re-inserting, case (4)
contradicts property (2) for $S_0$, cases (6), (7) and (8) contradict properties
(3) and (4). Hence, we are left with case (1), and $s_x = s$ for all $x\in G$.
Again, we use induction to show $\phi(x) = \phi(e)\cdot x$, and get an isomorphism
\beq
G \,\ni\, g &\mapsto& (\phi_g: \; x\,\mapsto\, g\cdot x) \,\in\, \Isom_\S(G)
\eeq

{\bf (C)} It is easy to find a generating system $S_0$ of $G$ which fulfills all
properties of subtheorem (B), except for property (4):
$s_1^{s_2} = s_1$ is always true. We follow through the proof of subtheorem (B)
until case (8) cannot be contradicted. Assume it is realized, i.e.\ we find
$x\in G$, $s\in S_0$ with $\phi(xs) = \phi(x)\cdot s^{-1}$. Then, for each
$r\in S\setminus \{s, s^{-1}\}$ we must have $\phi(xsr) = \phi(x)\cdot
s^{-1}\cdot r^{-1}$, and from excluding all other cases and property (2) of
$S_0$ we further find $\phi(xs^2) = \phi(x)\cdot s^{-2}$. By induction and
using the fact that $S_0$ generates $G$, we show
\beq
\phi(s_1\,s_2\,\ldots\, s_n) &=& \phi(e)\cdot s_1^{-1}\, s_2^{-1}\,\ldots\, s_n^{-1},
\eeq
or, due to abelianness, $\phi(x) = \phi(e) \cdot x^{-1}$. Obviously, all these
bijections are indeed shared isometries:
\beq
d(\phi(x),\, \phi(y)) &=& d(x^{-1},\, y^{-1})
\;=\; ||x\,y^{-1}|| \qquad |\n{ abelianness}\\
&=& ||y^{-1}\,x|| \;=\; d(y,\, x) \qquad\qquad |\n{ $S_0$ is symmetric}\\
&=& d(x,\,y)
\eeq

Hence, we have
$\Isom_\S(G)$ isomorphic to $G\rtimes \fg{2}$ via
\beq
G\rtimes \fg{2} \,\ni\, (g,\, a) &\mapsto& (\phi_{(g,\,a)}: \; x\,\mapsto\, g\cdot
x^a) \,\in\, \Isom_\S(G)
\eeq

{\bf (D)} Once again, we follow through the proof of subtheorem (B). As $S_0$ is
minimal, property (1) is automatically fulfilled. And as each $s\in S_0$ has
order 2, the question $s_x = s$ or $s_x = s^{-1}$ is trivial, as $s^{-1} = s$.
Hence, we get the usual isomorphism
\beq
G \,\ni\, g &\mapsto& (\phi_g: \; x\,\mapsto\, g\cdot x) \,\in\, \Isom_\S(G).
\eeq
\end{Proof}

From now on, we will restrict to the symmetric case $\S = \S_\n{sym}$.

\begin{Example}
For groups $G$ with central elements it can be difficult to find a generating
system $S_0$ satisfying the properties of Theorem \ref{THE_bartholdi}.B, but
typically it is still possible. Take for example:
\beq
G &=& \left\langle a,b,c\;|\;[a,c],[b,c]\right\rangle \;\cong\; (\Z \ast
\Z)\times \Z\\
S_0 &=& \left\{a^{\pm 1}, (bc)^{\pm 1}, (ab)^{\pm 1}\right\}.
\eeq
\end{Example}

\begin{Example}
The same accounts for groups with 2-torsion. For example, it is easy to
calculate by hand
\beq
\Isom_{\S}(\fg{2}) &\cong& \fg{2},
\eeq
just as Theorem \ref{THE_bartholdi}.D mentions; but not $\fg{2} \rtimes \fg{2}$, as one might
think from Theorem \ref{THE_bartholdi}.C. Indeed, as inversion is the trivial operation in each
group of exponent 2, we have $\Isom_{\S}(\fg{2})^n \cong (\fg{2})^n$ in the
abelian case, contrary to Theorem \ref{THE_bartholdi}.C.
\end{Example}

\begin{Example} \label{EXA_inverting_group_1}
For groups of the form $G \,=\, G_0\,\rtimes\, \fg{2}$ with $\fg{2}$ acting via
inversion (written multiplicatively) on a f.g.\ group $G_0$ (which subsequently
must be abelian), each element $(g,\,-1)$ with $g\,\in\, G_0$ has torsion $2$.
Given a minimal generating system $S_0$ of $G_0$, we can use
\beq
S &:=& \{(g,\,-1):\; g \,\in\, S_0\} \,\cup\, \{(e,\,-1)\}
\eeq
to apply Theorem \ref{THE_bartholdi}.D. And, just as it states, the inversion is
not a shared isometry in this case: Let $G_0$ be any f.g.\ group with at least
one element $s\,\in\, G_0$ with $s^2 \neq e$, $S_0$ a finite generating system of
$G_0$ with $s\,\in\, S_0$, and $S' \,:=\, S_0 \,\cup\, \{(s,\,-1)\}$, which
generates $G \,= \,G_0\rtimes \fg{2}$. Then holds $d\big((e,\,-1),\,
(s,\,1)\big) \,=\, 1$, as $(e,\,-1)\cdot (s,\,-1) \,=\, (s,\, 1)$, but 
\beq
d\left((e,\,-1)^{-1},\, (s,\,1)^{-1}\right) \;\;=\;\; d\left((e,\,-1),\, (s^{-1},\,
1)\right) \;\; > \;\; 1,
\eeq
because $(s^{-1},\, -1)\, \notin\, S'$. ($s \,\neq\, s^{-1}$, and $(s,\,-1)^{-1} \,=\, (s,\,-1)$.)
\end{Example}

\begin{Example} \label{EXA_inverting_group_2}
Similar to Example \ref{EXA_inverting_group_1}, consider a group
$G\,=\,G_0\,\rtimes\, H$, where a f.g.\ group $H$ acts on the f.g.\ abelian
group $G_0$. The action shall be given by a non-trivial homomorphism
$\alpha:\,H\,\rightarrow\,\fg{2}$, where $\fg{2}$ acts on $G_0$ by inversion. 
Furthermore, let $S_0$ be an arbitrary finite generating system of $G_0$, and
let $S_H$ be a finite generating system for $H$, such that there are no two
elements $s,\,t\,\in\,S$ with $s\,\neq\,t^{\pm 1}$ and $s^t\,=\,s^{\pm 1}$ or
$s^2\,t^2 = e$. Finally, let $h_0\,\in\,S_H$ be an element with $h_0^4\,\neq\,
e$. Then we can define a finite generating system
\beq
S_0 &:=& S_H \;\cup\; \big\{g h_0\;:\;g\,\in\,S_0\big\}
\eeq
from which we choose a minimal subsystem $S\,\subset\,S_0$. Some simple
calculations then show that the generating system $S$ fulfills the
requirements for Theorem \ref{THE_bartholdi}.B, and we conclude:
\beq
\Isom_\S(G_0\,\rtimes\,H) &\cong& G_0\,\rtimes\, H
\eeq
In particular, this accounts for the group
\beq
\Z\;\rtimes\;\Z\;\;=\;\;\langle x,\, y\;:\; x^y\,=\,x^{-1}\rangle &\cong&
\langle y,\,z\;:\; y^2 \,=\,z^2\rangle.
\eeq
\end{Example}

Considering the proof of Theorem \ref{THE_bartholdi} and the above examples, we
are confident that the following statements can be proven just by application of
more arduous combinatorics:

\begin{quote}
(A) Let $G$ be a f.g.\ group, and let $\S$ be the family of all symmetric
generating systems of $G$. Then $\Isom_\S(G) \cong G \rtimes \fg{2}$ if and only
if $G$ is non-trivial, abelian, and not of exponent 2; $\Isom_\S(G) \cong G$
otherwise.

(B) The shared Clifford isometries (i.e.\ those shared isometries $\phi$ with
constant $d(x,\,\phi(x))$ for all $x\in G$) always constitute a
group, which is isomorphic to $G$.
\end{quote}

\begin{Lemma} \label{LEM_fine_transfer}
Let $G$, $H$ be f.g.\ groups, $\S_G$, $\S_H$ families of generating systems of
$G$, $H$. If there is a bijection $\eta: G\rightarrow H$ such that
\begin{itemize}
\item for each $S_G\,\in\, \S_G$ there is $S_H\,\in\, \S_H$ which makes $\eta:
\Cay(G,\, S_G) \rightarrow \Cay(H,\, S_H)$ an isometry, and
\item for each $S_H\,\in\, \S_H$ there is $S_G\,\in\, \S_G$ which makes $\eta^{-1}:
\Cay(H,\, S_H) \rightarrow \Cay(G,\, S_G)$ an isometry.
\end{itemize}
Then $\Isom_{\S_G}(G)$ and $\Isom_{\S_H}(H)$ are isomorphic.

In particular, in the situations of Theorem \ref{THE_bartholdi}.A, B, or D, or
when $G$ and $H$ are both f.g.\ abelian without 2-torsion (case (C)), 
then $G$ and $H$ are isomorphic.
\end{Lemma}
\begin{Proof}
Define
\beq
\eta^*:\; \Isom_{\S_G}(G) &\rightarrow& \Isom_{\S_H}(H)\\
\phi &\mapsto& \eta\circ \phi \circ \eta^{-1}.
\eeq
This is well-defined: For each $S_H\in \S_H$ choose $S_G\in \S_G$ such that
$\eta$ is an isometry. Then $\eta\circ \phi\circ \eta^{-1}: H\rightarrow H$ is
an isometry as well---vice versa for $(\eta^*)^{-1} \,:=\, \eta^{-1}\circ
\cdot\circ \eta$. Hence, $\eta^*$ is a bijection, and, as one easily computes,
indeed an isomorphism between groups.

In the cases (A), (B) and (D), we may directly conclude $G\,\cong\, H$. In the abelian
case we just have $G\rtimes\fg{2} \,\cong\, H\rtimes\fg{2}$, but, as $G$ and $H$ are
without 2-torsion, $G$ and $H$ must be isomorphic as well.
\end{Proof}

\begin{Example}
Let us take a look at the three commensurable groups $G_1 \,=\, \Z$, $G_2 \,=\,
\Z\rtimes \fg{2}\,\cong\,\fg{2}\ast\fg{2}$, and $G_3 \,=\, \Z\times\fg{2}$. For
$G_1$, choose $S_0 = \{\pm 1,\, \pm 2\}$, and apply Theorem
\ref{THE_bartholdi}.C; for $G_2$ use Theorem \ref{THE_bartholdi}.D
(c.f. previous example); for $G_3$ apply a direct calculation\footnote{In this
case it suffices to find the isometries for the standard generating set, the
Cayley graph of which is a ladder. The cardinality of the second neighborhood of
an edge in this graph depends on the order of its generating element, but must
be preserved under isometries. This allows for a simple case distinction.}.
Then we find
\beq
\Isom_\S(G_1) &\cong& \Z\rtimes \fg{2} \\
\Isom_\S(G_2) &\cong& \Z\rtimes \fg{2} \\
\Isom_\S(G_3) &\cong& (\Z \rtimes \fg{2}) \times \fg{2} \;\;\cong\;\;
G_3\rtimes \fg{2}.
\eeq
We note that the resulting shared-isometry groups can be isomorphic, but might
as well be just commensurable. And, as $G_1$ and $G_2$ are not isomorphic,
we note that there cannot be a bijection $\eta: G_1 \rightarrow G_2$ as in Lemma
\ref{LEM_fine_transfer}. The canonical inclusion $i: G_1 \hookrightarrow G_2$
however might provide a deeper insight - it is a rough isometry for
several generating systems.
\end{Example}

\section{Shared Rough and Quasi-Isometries} \label{SEC_shared_rough}

Sometimes it is possible to directly translate a proof into the rough context.
This will be our goal for this section: To roughificate the proof of Theorem
\ref{THE_bartholdi}.

\begin{Definition} \label{DEF_optimal}
Let $G$ be a f.g.\ group. We call a family $\S$ of finite generating
systems of $G$ {\em optimal} if $\URIsom_\S(G) \,\cong\, G$, and {\em
quasi-optimal} if $\UQIsom_\S(G)\,\cong\, G$.
\end{Definition}

\noindent
Note that quasi-optimality is the stronger of both notions, because
\beq
\URIsom_\S(G)\,\subset\,\UQIsom_\S(G).
\eeq
Each translation from the left with an
element of $G$ is a shared isometry, and hence we have
\beq
G \;\;\leq\;\; \Isom_\S(G) \;\;\subset\;\; \URIsom_\S(G) \;\;\subset\;\; \UQIsom_\S(G).
\eeq
If $\S$ is optimal, we also find $\Iden_\S(G)$ to be trivial.

\begin{Lemma}[Optimality Lemma] \label{LEM_optimal}
Let $G$ be a finitely generated group with $x^y \,\neq\, x^{-1}$ for all
$x,\,y\,\in\, G$, unless $x\,=\,x^{-1}$. Let $G$ be non-abelian, or of exponent
2. Let $\S$ be a family of finite generating systems of $G$ with the following
{\bf Property \ref{LEM_optimal}}:
\begin{itemize}
\item
For each $g,\, h \,\in\, G$ with $g\,\neq\, h^{\pm 1}$ and each $R\,\in\,\N^*$
there is \\ $S \,=\, S(g,\,h,\,R) \,\in\, \S$ such that $g\,\in\, S$ and
$||h||_S\,\geq\, R$, or vice versa.
\end{itemize}
Then $\S$ is quasi-optimal (and thus optimal).
\end{Lemma}
\begin{Proof}
Let $\lambda,\,\epsilon \,\geq\, 0$, $\phi\, \in\,
(\lambda,\,\epsilon)\eIsom_\S(G)$, and $x,\, y\,\in\, G$ be arbitrary, let $z
\,:=\, y^{-1}\cdot x$ and define 
\beq
z' &:=& \phi(y)^{-1}\cdot \phi(x) \qquad \Rightarrow \qquad ||z'||_S \;=\;
d_S(\phi(y),\, \phi(x)).
\eeq
for all $S\,\in\, \S$. Now assume $z'\,\neq\, z^{\pm 1}$. Then there is
$S \,=\, S\big(z, \,z',\, (1\,+\,\epsilon)\cdot (1~\!+~\!\lambda)\big)$,
such that $\phi$ is still a $(\lambda,\,\epsilon)$-quasi-isometry, and it holds
\begin{itemize}
\item
either $z\,\in\, S$, then $||z'|| \,=\, d_S(\phi(y),\, \phi(x)) \,\leq\,
\lambda\,d_S(y,\,x) \,+\, \epsilon \,=\, \lambda \,+\, \epsilon$, but $||z'||_S
\,>\, \lambda\,+\,\epsilon$: contradiction,
\item
or $z'\,\in\, S$, then $d_S(\phi(y),\, \phi(x)) \,=\, 1$, hence $d_S(y,\, x)
\,\leq\, \lambda \,+\,\lambda\,\epsilon$, but $||z||_S \,>\, \lambda \,+\,
\lambda\,\epsilon$: contradiction!
\end{itemize}
Thus, $\phi(x) \,=\, \phi(y) \cdot (y^{-1}\cdot x)^{\pm 1}$, or (after
substitution): $\phi(yx) \,=\, \phi(y)\cdot x^{\pm 1}$. The sign might still
depend on $x$ and $y$, which we exclude in the next step.

Let $c \,:=\, \phi(e)$, and assume there are $x,\, y\,\in\, G$ with $\phi(x)
\,=\, c\,x \,\neq\, c\,x^{-1}$, but $\phi(xy) \,=\, c\, (xy)^{-1}\,\neq\,c\,x\,y$. Then
\beq
c\cdot (xy)^{-1} \;\;=\;\; \phi(xy) \;\;=\;\; \phi(x)\,y^a \;\;=\;\; c\,x\,y^\alpha
\eeq
for some $\alpha\,=\,\pm 1$, hence $x\,y^\alpha \,=\, y^{-1}\,x^{-1}$.
If $\alpha\,=\,+1$, we have $(xy)^2 \,=\, e$, and hence $\phi(xy) \,=\, c\,
(xy)$. If $\alpha\,=\,-1$, we have $x^y \,=\, x^{-1}$, which contradicts our
premise, unless $x\,=\,x^{-1}$. However, if $x\,=\,x^{-1}$, we have
$\phi(x)\,=\,c\,x^{-1}$.

We conclude that $\phi(x) \,=\, c\,x$ for all $x\,\in\, G$, or $\phi(x) \,=\,
c\, x^{-1}$ for all $x\,\in\, G$. The latter case leads to
\beq
c\,y\, x \;\;=\;\; \phi\big(x^{-1}y^{-1}\big) \;\;=\;\; \phi\big(x^{-1}\big)\,
y^\beta \;\;=\;\; c\,x\,y^\beta
\eeq
for all $x,y\in G$, and some $\beta\,=\,\pm 1$. Again, the case $\beta\,=\,-1$
leads to $y^x\,=\,y^{-1}$, which we excluded, unless $y\,=\,y^{-1}$. So both
cases for $\beta$ lead to the conclusion that $G$ must be abelian. Indeed, in
the abelian case, the inversion is a shared isometry of all symmetric finite
generating systems, and it is non-trivial if and only if $G$ is not of exponent
2.

Hence, $\UQIsom_\S(G)\, \cong\, \URIsom_\S(G)\,\cong\,G\rtimes \fg{2}$ if and
only if $G$ is abelian and not of exponent 2,
$\UQIsom_\S(G)\,\cong\,\URIsom_\S(G) \,\cong\, G$ otherwise.
\end{Proof}

\begin{Example} \label{EXA_torsion_vs_optimality}
No finite group has Property \ref{LEM_optimal}, as its diameter is limited.
Torsion in itself is an obstruction to it: Let $G$ have Property
\ref{LEM_optimal}, then each element of $G$ is either torsionsfree, or of
exponent 1, 2, 3, 4, or 6 -- these are those arguments for which the Euler
totient function $\varphi$ is 2 or less (\cite{Sloane_totient}): Let $x\,\in\,G$
be an element with $x^n \,=\, e$. If $\varphi(n) \,>\, 2$, we can choose two
different generators $a,\,b$ of $\fg{n}$, and hence $x^a$ and $x^b$ are powers
of each other, and yet $(x^a) \,\neq\, (x^b)^{\pm 1}$. Still, there might be
other optimal or quasi-optimal families for these groups.
\end{Example}

\begin{Example} \label{EXA_optimal_free}
Let $F_n$ be the free group generated by $S_0$ with $\#S_0\,=\,n\,\geq\,2$. Let
$g,\,h\,\in\,F_n$, $g\,\neq\,h^{\pm 1}$, and $R\,\in\,\N^*$ be arbitrary. Assume
$h$ is not a power of $g$ and not neutral (otherwise switch them; both cannot
happen as $F_n$ is torsionfree). If $g\,=\,e$, choose $x\,\in\,S_0$ such that $h$
is not a power of $x$, otherwise let $x\,=\,g$. Let $P$ be the maximum
of $R$ and the wordlength of $h$ in $S_0$. Define
\beq
\!\!\!\!S(g,\,h,\,R) &:=&
\{x\}\;\cup\;\left\{x^{\,(P+1)^j}\,s_j\;|\;s_j\,\in\,S_0\setminus\{x\},\,
j\,=\,1,\ldots \#(S_0\setminus\{x\})\right\}.
\eeq
The exponents $(P+1)^j$ are chosen such that any non-trivial product of the
elements $x^{\,(P+1)^j}\,s_j$ has large enough wordlength in $S_0$, that it
cannot equal $h$, at least for the first $R$ steps in the Cayley graph. After
this, the powers $x^{(P+1)^j}$ successively become available and ``free'' the
generators $s_j$ to generate each remaining element of $F_n$, such that
$||h||_S\,\geq\, R$. The family $\S$ of all these generating systems is
quasi-optimal due to Lemma \ref{LEM_optimal}.
\end{Example}

\begin{Example} \label{EXA_optimal_free_abelian}
In a similar way, we may define quasi-optimal generating families for free
abelian groups. We give the explicit example for $G\,=\,\Z$ (written
additively): Again, assume $h$ is not neutral and not a multiple of $g$. If $g$
is zero, let $P\,=\,1\,+\,(R\,\smax\,|h|)$, otherwise choose
$P\,\geq\,1\,+\,(R\,\smax\,|h|\,\smax\,|g|)$ and coprime to $g$. Then define
$S(g,\,h,\,R) \;:=\; \{g,\, P^2,\, P^3\,+\,1\}$.

For arbitrary f.g.\ free abelian groups, do this componentwise.
\end{Example}

The ``delayed generation method'' we applied in Examples \ref{EXA_optimal_free}
and \ref{EXA_optimal_free_abelian} can sometimes be generalized to other
f.g.\ groups: Choose a finite generating system $S_0$, then find a suitable
element $x\in G$ such that $x$ and $g$ together do not generate $h$. Add
$x^Ps_1$, $x^{P^2}s_2$, $x^{P^3}s_3$ and so on, after choosing $P$ large
enough and taking care for the group's relations: If e.g.\ holds
$x^Ps_1\,=\,s_3$, choose $P$ even larger, or change the sequence of the
generators.

\vspace*{11pt}

\begin{Question}
Is there a torsion-free group without Property \ref{LEM_optimal}, or which
does not admit an optimal generating family?
\end{Question}

\begin{Definition} \label{DEF_semishared}
Let $G$ and $H$ be f.g.\ groups, and let $\S_H$ be a family of generating
systems of $H$. We call a pair of maps $\eta:\,G\,\rightarrow\,H$ and
$\eta':\,H\,\rightarrow G$ an {\em $S_H$-semi-shared quasi-isometry} if there
are $\lambda,\,\epsilon\,\geq 0$ with:
\begin{itemize}
\item For each $S_H\,\in\,\S_H$ there is a finite generating system $S_G$ of
$G$ which makes $(\eta,\,\eta'):\,\Cay(G,\,S_G)\,\rightarrow\,\Cay(H,\,S_H)$ a
$(\lambda,\,\epsilon)$-quasi-isometry.
\end{itemize}
When we speak of an ``$\S_H$-semi-shared quasi-isometry
$\eta:\,G\,\rightarrow\, H$'' a suitable $\eta'$ shall always be implied.
\end{Definition}

\begin{Theorem} \label{THE_semishared_isomorph}
Let $G$ and $H$ be f.g.\ groups with $x^y \,\neq\,x^{-1}$ for all $x,\,y\,\in\,
G$ (resp. $H$), unless $x\,=\,x^{-1}$. Let $G$ and $H$ be  non-abelian, or of
exponent 2. Let $\S_G$ and $\S_H$ be quasi-optimal families of $G$ and $H$,
respectively, and let $\eta:\,G\rightarrow\,H$ be an $\S_H$-semi-shared
quasi-isometry, such that $\eta'$ is an $\S_G$-semi-shared quasi-isometry. Then
$G$ and $H$ are isomorphic as groups.
\end{Theorem}
\begin{Proof}
We first note that $\eta\circ\eta':\, H\rightarrow H$ is an $\S_H$-uniform
quasi-isometry, and hence it is given by multiplication from the left with an
element $c\,\in\,H$. Consider $\eta'':\,G\rightarrow H$
given by $h \mapsto c^{-1}\cdot\eta(h)$. Then $(\eta'',\, \eta')$ is another
$\S_H$-semi-shared quasi-isometry, $(\eta',\,\eta'')$ is a $S_G$-semi-shared
quasi-isometry, and $\eta''\circ\eta'$ is the identity.

Conversely, $\eta'\circ\eta'':\,G\rightarrow G$ also is a multiplication from the
left with an element $c'\,\in\,G$. We easily find
\beq
(\eta'\,\eta''\,\eta'\,\eta'')(h) &=& (c')^2\cdot h\\
&=& (\eta'\cdot\id_H\cdot \eta'')(h) \;\;=\;\; c'\cdot h
\eeq
for all $h\,\in\,H$, thus $c'\,=\,e$, and consequently $\eta'\circ\eta''$ is the
identity as well.

Without loss of generality, and to ease our notation, we may assume that
$(\eta,\,\eta')$ already fulfills $\eta\circ \eta'\,=\,\id_H$ and
$\eta'\circ\eta\,=\,\id_G$.

Now define
\beq
\xi: \quad \UQIsom_{\S_H}(H) &\rightarrow& \UQIsom_{\S_G}(G)\\
\phi &\mapsto& \eta' \circ \phi \circ \eta, \quad \n{ and }\\
\xi': \quad \UQIsom_{\S_G}(G) &\rightarrow& \UQIsom_{\S_H}(H)\\
\psi &\mapsto& \eta \circ \psi \circ \eta'.
\eeq
$\xi$ is well-defined: For each $S_G\,\in\, \S_G$ choose $S_H\,\in\,\S_H$ such that
$(\eta,\, \eta')$ is a uniform quasi-isometry. Then for each $\phi\,\in\, \UQIsom_{\S_H}(H)$,
the map $\eta'\circ \phi\circ \eta : G\rightarrow G$ is a uniform quasi-isometry as well; the
same accounts for $\xi'$. Furthermore, we have
\beq
\xi(\phi_1) \circ \xi(\phi_2) \;\;=\;\; \eta'\,\phi_1\,\eta\,\eta'\,\phi_2\,\eta
&=& \xi(\phi_1\circ\phi_2) \\
\xi'(\psi_1)\circ\xi'(\psi_2) \;\;=\;\; \eta\,\psi_1\,\eta'\,\eta\,\psi_2\,\eta'
&=& \xi'(\psi_1\circ\psi_2) \\
\xi'(\xi(\phi)) \;\;=\;\; \eta\,\eta'\,\phi\,\eta\,\eta' &=& \phi \\
\xi(\xi'(\psi)) \;\;=\;\; \eta'\,\eta\,\psi\,\eta'\,\eta &=& \psi
\eeq
for all $g\,\in\, G$. This means that $(\xi,\,\xi')$ constitutes an isomorphism
between $\UQIsom_{\S_H}(H)\,\cong\,H$ and $\UQIsom_{\S_G}(G)\,\cong\,G$.
\end{Proof}

A similar theorem should hold in the abelian case.

We now want to weaken the hypothesis of Theorem \ref{THE_semishared_isomorph}. For this,
we will rework the proof of Lemma \ref{LEM_optimal}, which yields a
generalized form of homomorphism.

\begin{Lemma} \label{LEM_generalized_shared_isometry}
Let $G$ and $H$ be f.g.\ groups, and $\S_H$ a family of generating systems of $H$
satisfying Property \ref{LEM_optimal}. Let $\phi: \,G\,\rightarrow\, H$ be an
$\S_H$-semi-shared quasi-isometry with $\phi(e)\,=\,e$. Then $\phi$ fulfills
$\phi(g h)\,=\,\phi(g)\cdot\phi(h)^{\pm 1}$ for all $g,\,h\,\in\,G$, where the
sign might depend on $g$ and $h$.
\end{Lemma}
\begin{Proof}
Let $x,\,y\,\in\,G$ be arbitrary, $z\,=\,y^{-1}\,x$, and
$z'\,=\,\phi(y)^{-1}\,\phi(x)$. Assume $z'\,\neq\,\phi(z)^{\pm 1}$. Then
we may choose
$S_H\,=\,S\big(z',\,\phi(z),\,(\lambda^2\,+\,1)\cdot(\epsilon\,+\,1)\big)\,\in\S_H$
a suitable generating system to separate $z'$ from $\phi(z)$. Then we have
\beq
\big|\big|z'\big|\big|_{S_H} &=& d_{S_H}(\phi(y),\, \phi(x))\\
&\leq& \lambda\cdot d_{S_G}(y,\,x)\;+\;\epsilon\;\;=\;\; \lambda\cdot d_{S_G}(e,\,z)\;+\;\epsilon\\
&\leq& \lambda^2\cdot
d_{S_H}(\phi(e),\,\phi(z))\;+\;\lambda^2\,\epsilon\;+\;\epsilon \\
&=& \lambda^2\cdot \big|\big|\phi(z)\big|\big|_{S_H} +\;\lambda^2\,\epsilon\;+\;\epsilon,
\eeq
and, similarly:
\beq
\big|\big|\phi(z)\big|\big|_{S_H} &\leq&
\lambda^2\cdot \big|\big|z'\big|\big|_{S_H} +\;\lambda^2\,\epsilon\;+\;\epsilon,
\eeq
Now one of $||z'||_{S_H}$ and $||\phi(z)||_{S_H}$ is $1$, while the other is
larger than $\lambda^2\,+\,\lambda^2\,\epsilon\,+\,\epsilon$, contradiction. Hence,
$z'$ is $\phi(z)^\alpha$ for some suitable $\alpha\,=\,\pm 1$, which depends on
$x$ and $y$. Substituting $y\,=\,g$ and $z\,=\,h$ yields
$\phi(g h)\,=\,\phi(g)\cdot\phi(h)^{\pm 1}$.
\end{Proof}

Note that it is always possible to switch from an arbitrary semi-shared
quasi-isometry $\phi$ to one with $\phi(e)\,=\,e$ by a simple translation. The
translation even preserves the constants $\lambda$ and $\epsilon$ of the
quasi-isometry.

\begin{Theorem} \label{THE_semishared_commensurable}
Let $G$, $H$, and $\phi:\,G\,\rightarrow\, H$ be as in Lemma
\ref{LEM_generalized_shared_isometry}. \\
Assume one of the following statements holds:
\begin{enumerate}
\item $G$ admits a generating system $S$ such that $\big(\phi(s)\big)^2\,=\,e$
for each $s\,\in\,S$.
\item $G$ admits a generating system $S$ such that:
\begin{enumerate}
\item There is no $x\,\in\,\phi(S\,\cup\,S^{-1})$, with $x^2\,=\,e$.
\item There are no $x,\,y\,\in\,\phi(S\,\cup\,S^{-1})$, $x\,\neq\,
y^{\pm 1}$, with $x^2\,=\,y^2$.
\item There are no $x,\,y\,\in\,\phi(S\,\cup\,S^{-1})$, $x\,\neq\,
y^{\pm 1}$, with $(xy)^2\,=\,e$.
\item There are no $x,\,y\,\in\,\phi(S\,\cup\,S^{-1})$, $x\,\neq\,
y^{\pm 1}$, with $x^y\,=\,x$.
\item There are no $x,\,y\,\in\,\phi(S\,\cup\,S^{-1})$, $x\,\neq\,
y^{\pm 1}$, with $x^y\,=\,x^{-1}$.
\item There are at least two distinct elements in $S$, which are not
inverses of each other.
\end{enumerate}
(In particular, $G$ is not abelian.)
\end{enumerate}
Then $G$ and $H$ are commensurable up to finite kernels
(Definition \ref{DEF_commensurable}).
\end{Theorem}

\goodbreak

\begin{table}
\center
\begin{tabular}{|r|c|c|c|c|l|l|}
\hline
Nr&$\alpha$&$\beta$&$\gamma$&$\delta$& $x^\alpha\cdot y^\beta \;=\; (x\cdot
y^\gamma)^\delta$ & contradiction \\ \hline
1&$+$&$+$&$+$&$+$&  --                  & {\bf no} \\
2&$+$&$+$&$+$&$-$&  $(xy)^2\;=\;e$      & yes (c) \\
3&$+$&$+$&$-$&$+$&  $y^2\;=\;e$         & yes (a) \\
4&$+$&$+$&$-$&$-$&  $x^y\;=\;x^{-1}$    & yes (e) \\ 
5&$+$&$-$&$+$&$+$&  $y^2\;=\;e$         & yes (a) \\
6&$+$&$-$&$+$&$-$&  $x^y\;=\;x^{-1}$    & yes (e) \\ 
7&$+$&$-$&$-$&$+$&  --                  & {\bf no} \\
8&$+$&$-$&$-$&$-$&  $(xy^{-1})^2\;=\;e$ & yes (c) \\
9&$-$&$+$&$+$&$+$&  $x^2\;=\;e$         & yes (a) \\
10&$-$&$+$&$+$&$-$& $y^x\;=\;y^{-1}$    & yes (e) \\ 
11&$-$&$+$&$-$&$+$& $x^2\;=\;y^2$       & yes (b) \\
12&$-$&$+$&$-$&$-$& $x^y\;=\;x$         & yes (d) \\ 
13&$-$&$-$&$+$&$+$& $x^{-2}\;=\;y^2$    & yes (b) \\
14&$-$&$-$&$+$&$-$& $x^y\;=\;x$         & yes (d) \\ 
15&$-$&$-$&$-$&$+$& $x^2\;=\;e$         & yes (a) \\
16&$-$&$-$&$-$&$-$& $y^x\;=\;y^{-1}$    & yes (e) \\
\hline
\end{tabular}
\caption{The sixteen cases of the proof of Theorem
\ref{THE_semishared_commensurable}.3. For convenience, we use $x\,=\,\phi(s)$
and $y\,=\,\phi(t)$.}
\label{TAB_sixteen_cases_commensurability_theorem}
\end{table}

\noindent
\begin{Proof}
Due to Lemma \ref{LEM_generalized_shared_isometry} we have in each case
\beq
\phi(g\,h)&=&\phi(g)\cdot\phi(h)^{\sigma(g,\,h)}
\eeq
with $\sigma(g,\,h)\,\in\,\{\pm 1\}$ for any $g\,\in\,G$ and $h\,\in\,H$.
Observe that $\sigma(g,\, e)\,=\,\sigma(e,\,g)$ $\,=\,+1$.
If $\phi(h)$ is neutral or of order 2, we choose $\sigma(g,\,h)$ to be $+1$
without loss of generality.
We next show that under both hypothesis $\phi$ must be a homomorphism.
Due to Proposition \ref{PRO_homomorphism_and_qi} $G$ and $H$ then must be
commensurable up to finite kernels.

{\bf (1)} We trivially have
\beq
\phi(g\,s)&=&\phi(g)\cdot\phi(s)
\eeq
for any $g\,\in\,G$ and $s\,\in\,S$. By induction, $\phi$ must be a homomorphism.

{\bf (2)} Let $g\,\in\,G$ and $s,\,t\,\in\,S$ be arbitrary, $s\,\neq\,t^{\pm 1}$.
We make use of the associative law:
\beq
\phi(g\,s\,t) &=& \phi(g)\,\cdot\,\phi(s)^\alpha\,\cdot\,\phi(t)^\beta\\
&=& \phi(g)\,\cdot\,\big(\phi(s)\,\cdot\,\phi(t)^\gamma\big)^\delta
\eeq
for some $\alpha,\,\beta,\,\gamma,\,\delta\,=\,\pm 1$. The sixteen possible
cases resolve as in Table \ref{TAB_sixteen_cases_commensurability_theorem}.
Fourteen cases subsequently contradict our premise.
Both remaining cases 1 and 7 demand $\alpha\,=\,\sigma(g,\,s)\,=\,+1$, for all
$g\,\in\,G$ and $s\,\in\,S$, so we have
\beq
\phi(g\,s) &=& \phi(g)\cdot \phi(s),
\eeq
and, again by induction, $\phi$ must be a homomorphism.
\end{Proof}

\begin{Corollary} \label{COR_semishared_quotient}
Let $G$, $H$, and $\phi$ be as in Theorem \ref{THE_semishared_commensurable},
and let $H$ be non-abelian, or of exponent 2. In addition, $x^y\,\neq\,x^{-1}$
shall hold for all $x,\,y\,\in\,H$ with $x\,\neq\,x^{-1}$. Then $H$ is the
quotient of $G$ by the finite subgroup $\ker\phi\,\trianglelefteq\,G$.
\end{Corollary}
\begin{Proof}
Lemma \ref{LEM_optimal} ensures that $\phi\circ\phi':\,H\rightarrow H$ is
given by multiplication with a fixed element of $H$, and in particular, $\phi$
must be surjective. From the proof of Theorem \ref{THE_semishared_commensurable}
we know that $\phi$ is a homomorphism with finite kernel. Using the First
Isomorphism Theorem (\cite{Bosch}, Korollar 1.2.7), we see
$H\,=\,\im\phi\,\cong\,G/\ker\phi$.
\end{Proof}

We present a very simple example to demonstrate that Theorem
\ref{THE_semishared_commensurable} is not empty.

{ 
\begin{Example}
Choose $H\,=\,\Z$ and $G\,=\,H\,\rtimes\,C_4$ with inversion as
action, i.e.
\beq
(x,\,0)^{(0,\,0,\,t)} &:=& \left\{
\begin{array}{r}
(x,\,0)\quad \n{if} \quad t\,\equiv\,0,\,2\;\;(\n{mod } 4)\\
(x^{-1},\,0)\quad \n{if} \quad t\,\equiv\,1,\,3\;\;(\n{mod } 4)
\end{array} \right.
\eeq
for any $x\,\in\,\Z$ and $t\,\in\,C_4$.
Note that the
restriction in torsion for Property \ref{LEM_optimal} (see Example
\ref{EXA_torsion_vs_optimality}) only applies to $H$, not to $G$.
Choose $\S_H$ like in Example \ref{EXA_optimal_free_abelian} and the following
generating system for $G$:
\beq
S &:=& \{(1,\,1),\,(2,\,1)\}~.
\eeq
$S$ fulfills the requirements of case (2) in Theorem
\ref{THE_semishared_commensurable}, as one easily calculates.
Finally, our choice for the quasi-isometry is the obvious one:
\beq
\phi: \qquad G&\twoheadrightarrow & H\\
(g,\,t) &\mapsto& g\\
\n{and}\qquad \phi': \qquad H&\hookrightarrow & G\\
g &\mapsto& (g,\,0)~.
\eeq
We are only left to show that $(\phi,\,\phi')$ is an $\S_H$-semi-shared
quasi-isometry.

Set $\lambda\,=\,1$ and $\epsilon\,=\,1$. Given any generating set
$S_H\,\in\,\S_H$ we have to show that there is a finite generating set $S_G$ of
$G$ which makes $(\phi,\,\phi')$ a $(\lambda,\,\epsilon)$-quasi-isometry.
Assume
\beq
S_H &=& \{g,\;P^2,\;P^3\,+\,1\}
\eeq
with $g,\,P$ as in Example \ref{EXA_optimal_free_abelian}. Choose
\beq
S_G &:=& \big\{\phi'(g),\;\phi'(P^2),\;\phi'(P^3\,+\,1),\;(0,\,1),\;(0,\,2)\big\}~.
\eeq
Adding the inversion as action doesn't change the metric on the subgroup
$(\cdot,\,0)\,\subset\,G$, because all inverses of the generating elements
$\phi'(g)$, $\phi'(P^2)$ and $\phi'(P^3\,+\,1)$ are already included in the
system $S_G\cup S_G^{-1}$, and because $\phi'$ is a homomorphism, such that
these inverses are included in $S_H\cup S_H^{-1}$ as well. All remaining
elements can be reached within one step.
\end{Example}
}

Unfortunately, Proposition \ref{PRO_normal_subgroup} is not yet strong enough
to constitute a reversal of Corollary \ref{COR_semishared_quotient}. Still, we
are confident to find a sustainable connection between semi-shared
quasi-isometries and quotients of finite kernel. Through the means of
residual finiteness, it might then be possible to finally find a 
perfectly fitting geometrical equivalence relation which equals
commensurability.

Georg-August-Universit\"at G\"ottingen, Germany\newline
eMail \verb|lochmann@uni-math.gwdg.de|

\end{document}